# Dynamic Simulation Method for Low-permeability Reservoirs with Fracturing-flooding Based on a Dual-Porous and Dual-Permeable Media Model

Xiang Wang[a], *, Wenjie Yu [a], Yixin Xie [b], Yanfeng He [a], Hui Xu [a], Xianxiang Chu [a], Changfu Li [a]

[a] *School of Petroleum and Natural Gas Engineering, Changzhou University, Changzhou 213164, China*

[b] *School of Mechanical Engineering, Changzhou Institute of Technology, Changzhou 213032, China*

**Abstract:**

The fracturing-flooding technology is a new process for the development of low-permeability oil reservoirs, achieving a series of successful applications in oilfield production. However, existing numerical simulation methods for pressure drive struggle to efficiently and accurately simulate the dynamic changes in reservoir properties during the fracturing-flooding process, particularly the expansion and closure of fractures within the reservoir. This paper introduces a Darcy flow model with dual-porous and dual-permeable characteristics based on seepage mechanics theory, utilizing two sets of rock stress-sensitive parameter tables to describe the physical property changes of the matrix and fractures during the fracturing-flooding process. Different parameters are set for the X and Y directions to characterize the anisotropic features of the reservoir. A numerical simulation method aimed at dynamic analysis of fracturing-flooding is established, along with an automatic history fitting method based on the CMA-ES algorithm to derive rock mechanics parameters that align with actual block conditions. Verification shows that the fitted model can accurately describe the dynamic changes of fractures within the reservoir. Compared to existing numerical simulation methods for fracturing-flooding, the modeling process is simpler and more efficient. Simulations of different timing for transitioning to conventional water injection indicate that the optimal development effect can be achieved when the reservoir pressure coefficient drops to around 1.2.

**Keywords:**

fracturing-flooding; numerical simulation; dual-porous and dual-permeable media; automatic history matching

## Introduction

Fracturing-flooding is a new technology for developing low-permeability oil reservoirs. It utilizes high-pressure water injection equipment to inject a large volume of water into the reservoir at high pressure within a short period, quickly replenishing the formation energy. This is followed by a period of no injection and no production, allowing the pressure to diffuse fully [1]. By creating a high-pressure zone near the injection well, this method activates existing closed fractures or generates micro-fractures within the reservoir, causing these fractures to extend and propagate[2]. This results in a rapid replenishment of reservoir energy, an increase in reservoir pressure, and subsequently enhances the production pressure differential and output volume of the oil well. In recent years, major oilfields in China, such as Shengli [3], Changqing [4], and Daqing [5], have begun applying fracturing-flooding technology on-site, achieving excellent results. The injection capacity of water injection wells in the test blocks has significantly improved, bottom-hole pressures of production wells have increased, and both liquid and oil production volumes have shown obvious improvement.

Unlike conventional waterflooding and traditional fracturing, fracturing-flooding is a process that couples fracture extension with reservoir seepage in real time. Establishing a convenient and efficient numerical simulation method for fracturing-flooding is of great significance for clarifying the mechanisms by which fracturing-flooding improves recovery rates, exploring the boundaries of technical policies, and further enhancing the effectiveness of fracturing-flooding. There is limited specialized research on fracturing-flooding models. Lu established a numerical model for two-phase flow of oil and water in the matrix and fracture media to investigate engineering issues such as rapid water cut increase and low water drive efficiency during the fracturing-flooding process after waterflooding in low-permeability sandstone reservoirs. He proposed an optimization method for fracture parameters in low-permeability sandstone reservoirs under fracturing-flooding conditions. The results indicated that the main controlling factors affecting rock fracture complexity include quartz content, Young's modulus, and peak stress. When fracturing production wells in the direction of the maximum horizontal principal stress, it is recommended to slow down the advance speed of the waterflooding front and reduce the risk of water channeling in injection-production wells [6]. Wang and colleagues considered flow rate changes caused by fracture closure based on Duhamel's principle, simulated fluid flow in variable mass radial flow fractures using the point source method and obtained the pressure response of injection wells inducing horizontal fractures. They noted that, compared to conventional fracturing wells, fracturing-flooding wells, which lack proppant reinforcement, exhibit changing fracture characteristics over time during shut-in tests [7]. Zhang and colleagues proposed a simulation model for matrix permeability variations that considers the mechanism of fracturing-flooding. Using numerical simulation software, they analyzed the effects of different parameters on development outcomes and established two types of reservoir numerical models: a pure matrix system and a fracture-matrix system, both considering the mechanism of fracturing-flooding [8]. Wang et al. established a multi-process multiphase flow model for fracturing, soaking, and production, considering capillary pressure, permeability pressure, membrane effects, and elastic properties. They inverted fracture parameters based on field fracturing pressure data and simulated the dynamic changes in the reservoir during fracturing, soaking, and production stages. The study showed that as soaking time increased, the cumulative production increment first rose rapidly and then approached a stable value. The inflection points of this change corresponded to the optimal soaking time [9]. Fan et al. proposed a multi-cycle expansion-recompaction geomechanical model that uses quantitative simulation techniques to describe the periodic deformation of rocks during steam injection processes and their geomechanical behavior. This model quantitatively characterizes the dynamic changes in porosity as reservoir pressure varies cyclically [3]. Shi developed typical reservoir models for different fault-block oil reservoirs and studied the distribution patterns of remaining oil during high-water cut periods. Based on summarizing the characteristics of remaining oil distribution, the study evaluated the applicability of fracturing-flooding techniques for three types of reservoirs: fractured, porous, and composite. The research found that in fractured reservoirs, the fracture system is the main flow channel, with higher affected volume and utilization compared to the matrix system. In porous reservoirs, water cut gradually increases, and the injection front expands more noticeably in the vertical direction than in the horizontal plane. In composite reservoirs, the injection front of both matrix and fracture systems gradually expands, with a more pronounced vertical effect compared to the horizontal plane [10]. Cha computed the rock stress-strain relationship using the finite element method, performed polynomial projections through node network degrees of freedom, and constructed a virtual element function space using the degrees of freedom (node values) within and at the boundaries of the elements. The continuous damage model introduced a damage factor and combined it with stress-strain load parameters within the mesh to calculate and predict the damage behavior mechanism of fractures within the matrix[11]. Cui et al., based on a stress sensitivity mathematical model and existing research on the relationship between formation pressure and fracture permeability, developed a permeability model for fractures under fracturing-flooding conditions. They used finite difference methods to simulate fracturing-flooding experiments and field conditions to verify the model's accuracy. The study indicated that injection rate is positively correlated with fracture extension speed; under the same injection rate, fracture extension speed is faster near the wellbore and slower farther from the wellbore[12].

The above research and field production verification can be summarized to categorize existing fracturing-flooding simulation methods into three types. The first type involves artificially modifying the permeability of the near-well zone to simulate fractures. The problem with this method is that the designed model is static, meaning its properties do not change as the development plan progresses, making it difficult to describe the changes in reservoir properties during the fracturing-flooding process. The second type uses complex fracture network expansion simulation methods, which face the issue of excessive computational demands, making them unsuitable for practical production needs. The third type simulates fracture expansion through fracturing numerical simulations and then embeds this into the fracturing-flooding numerical model. The problem with this method is the clear distinction between fracturing and fracturing-flooding fractures; compared to fracturing, fracturing-flooding has a longer injection time, larger total water volume, and no proppant influence, making fracturing-flooding simulations more complex [13]. Therefore, this method also struggles to accurately describe the changes in fractures during the fracturing-flooding process.

This study aims to establish a numerical simulation method capable of simulating the dynamic expansion and closure of fractures during the fracturing-flooding process, while ensuring that the implementation process is relatively simple for rapid analysis of the technical policy boundaries of fracturing-flooding. To this end, this paper introduces a Darcy flow model with dual-porosity and dual-permeability characteristics based on seepage mechanics theory, employing two sets of rock stress sensitive parameter tables to describe the physical property changes of the matrix and fractures in the rock during the fracturing-flooding process. Different parameters are set for the X and Y directions to characterize the anisotropic features of the reservoir. An automatic history matching method based on the CMA-ES algorithm is designed, which takes the bottom hole pressure curves, injection rate curves, water cut curves, and fracture expansion ranges of each production well during the fracturing-flooding process as optimization objectives, while using reservoir physical property parameters, such as pore volume multipliers, as optimization variables. The CMA-ES algorithm is then invoked to solve for the optimal parameter combinations. Finally, the feasibility and accuracy of the model and the automatic history matching method are verified against an actual block, and the timing for transitioning to conventional water injection after fracturing-flooding in the example block is studied. The innovations of this paper are as follows: (1) A simplified method is adopted to describe the dual-porous and dual-permeable media and anisotropic characteristics in low-permeability reservoirs, significantly reducing the complexity of the model while ensuring its accuracy, making it more convenient and efficient for practical applications in fracturing-flooding numerical simulation; (2) The automatic history fitting method based on the CMA-ES algorithm significantly improves the accuracy of the rock mechanics parameter settings in the dual-porous and dual-permeable media model and reduces the workload of fitting.

This paper is divided into five sections. In the first section, the dynamic changes in reservoir properties during the fracturing-flooding process are detailed. The second section establishes a dual-porous and dual-permeable media model considering the anisotropic stress sensitive characteristics of the rock based on the laws of seepage mechanics and rock mechanics. The third section describes the design of the automatic history matching method based on the CMA-ES algorithm. In the fourth section, we validate the accuracy of the model through matching results from an actual block and then discuss the optimization design of subsequent fracturing-flooding schemes based on this model. Finally, we present some conclusions in the fifth section.

## 1 The dynamic process of fracturing-flooding

Fracturing-flooding is a process that couples fracture expansion with reservoir seepage in real-time. Generally, the changes in fracture morphology during the fracturing-flooding process can be divided into three stages[12]: when fracturing-flooding is implemented, the reservoir pressure exceeds the rock's fracture pressure, causing the fractures to initiate and rapidly extend; when the reservoir pressure is between the closure pressure and the fracture pressure, the fracture morphology remains unchanged; after production begins, the reservoir pressure continuously decreases, and when it falls below the closure pressure, the fractures will partially close. The fractures formed by fracturing-flooding can effectively provide energy to the reservoir, driving crude oil towards the wellhead and shortening the effective interaction distance between the oil and water wells [7]. As fracturing-flooding progresses, once the pressure in the reservoir rises to a certain inflection point, the permeability of the reservoir will experience a rapid increase, indicating that fractures have appeared within the reservoir after reaching this point. The larger the scale of the fractures, the more pronounced the pressure diffusion range, leading to a faster increase in pressure around the well and a greater production pressure differential. According to laboratory core experiments, as the pore pressure increases during the fracturing-flooding process, when the fracture pressure and damage pressure points are reached, the permeability of the reservoir significantly increases, with fracture widths expanding from tens of microns to hundreds of microns, as shown in Fig 1.1 [14]. When the net horizontal stress reaches -6 MPa, the fracture width is 56 μm; at -7 MPa, it is 185 μm; and at -7.5 MPa, it is 281 μm. Fracturing-flooding achieves the goal of weakening the adverse effects of reservoir stress sensitivity and opening reservoir fractures, thus significantly enhancing the seepage capacity of the reservoir by gradually increasing the injection pressure in the field. The higher displacement pressure during fracturing-flooding creates a fracture modification zone around the injection well, while areas farther away from the injection well do not generate fractures due to the pressure not reaching the fracture pressure point, remaining as the matrix zone [15].

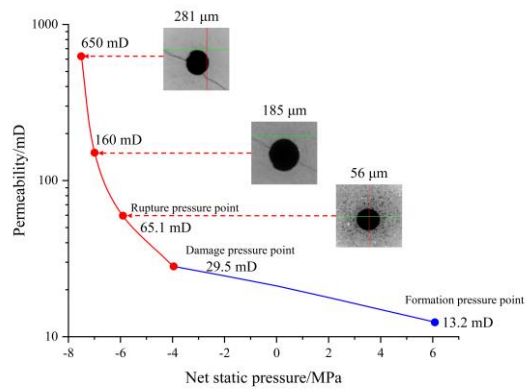

Fig 1.1 Net static pressure and permeability curve

Unlike conventional water injection and traditional fracturing, the dynamic changes in fractures during fracturing-flooding are more complex due to the rock fracturing and fracture expansion, influenced by factors such as long injection times, large total volumes, and the absence of proppants. The pressure changes and their impact range during the injection and production phases of fracturing-flooding are significant, requiring simultaneous consideration of the dynamic changes in physical properties in both the matrix and fracture zones. High-pressure water injection and oil-water displacement occur simultaneously and interact with each other, with the injection process triggering changes in geostress, modifying microfracture zones and their physical properties, which in turn affects the displacement effectiveness. The displacement process also reacts back on the geostress, thereby influencing the changes in fractures and reservoir properties [16].

## 2 Numerical simulation method for fracturing-flooding dynamics based on dual-porous and dual-permeable media model

As seen in the previous section, fracturing-flooding is a complex dynamic process, with its core being the description of the dynamic process of microfracture initiation, expansion, and recovery near the water well, as well as the dynamic interactions between fractures and the matrix. To address this, this study proposes a dynamic numerical simulation method for fracturing-flooding based on a dual-porous and dual-permeable media model. In this method, the dual-porous and dual-permeable media model is used to describe the interactions between the matrix and the fracturing-flooding fractures, while a rock stress sensitive model is employed to simulate the dynamic process of fracture opening and closure.

### 2.1 Numerical simulation basic model

The basic mathematical models involved in the oil-water two-phase seepage numerical model mainly include the motion equations, state equations, and continuity equations.

(1) Continuity Equation

The mass conservation equation, also known as the Euler continuity equation, is the mathematical expression of the law of mass conservation. The law of mass conservation in fluid mechanics indicates that the mass of fluid remains constant during motion, meaning that the difference between the mass of fluid flowing into a unit volume and the mass of fluid flowing out of the unit volume per unit time equals the increase in mass within the unit volume over that period.

$$-\frac{\partial}{\partial x}(\rho V_x) - \frac{\partial}{\partial y}(\rho V_y) - \frac{\partial}{\partial z}(\rho V_z) = \frac{\partial}{\partial t}(\rho \phi) \tag{2-1}$$

Where $\rho$ is the fluid density, kg/m$^3$; $\phi$ is the porosity; $V_x$, $V_y$, $V_z$ are the fluid velocities in the $x$, $y$, $z$ directions, m/s.

(2) Motion Equation

Darcy's equation reflects the fundamental laws of fluid motion in porous media, which serves as the motion equation in the reservoir numerical model.

The multiphase flow Darcy's law is expressed as:

$$Q_o = \frac{KK_{ro}A}{\mu_o L}\Delta p \tag{2-2}$$

$$Q_w = \frac{K K_{rw} A}{\mu_w L} \Delta p \qquad (2\text{-}3)$$

Where $Q_o$ is the oil phase flow rate, m³/d; $Q_w$ is the water phase flow rate, m³/d; $K$ is the air permeability, mD; $A$ is the cross-sectional area, m²; $L$ is the length, m; $\Delta p$ is the pressure difference, MPa; $K_{ro}$, $K_{rw}$ are the relative permeabilities of oil and water, respectively; $\mu_o$, $\mu_w$ are the viscosities of crude oil and water, respectively, MPa·s.

The motion equation follows Darcy's linear seepage law, and the form of Darcy's law motion equation is:

$$\begin{cases} V_x = \dfrac{-K_x}{\mu}\left(\dfrac{\partial p}{\partial x} - \rho g \dfrac{\partial D}{\partial x}\right) \\ V_y = \dfrac{-K_y}{\mu}\left(\dfrac{\partial p}{\partial y} - \rho g \dfrac{\partial D}{\partial y}\right) \\ V_z = \dfrac{-K_z}{\mu}\left(\dfrac{\partial p}{\partial z} - \rho g \dfrac{\partial D}{\partial z}\right) \end{cases} \qquad (2\text{-}4)$$

Where $V_x$, $V_y$, $V_z$ are the fluid velocities in the $x$, $y$, $z$ directions, respectively, m/s; $\rho$ is the fluid density, kg/m³; $g$ is the gravitational acceleration, m/s²; $D$ is the depth based on the vertical direction from the reference plane, m.

(3) State Equation

The state equation is a mathematical equation that describes how the state parameters of fluids and rocks change with pressure. It is usually represented by a compressibility coefficient. Before the oil and gas field is developed, the reservoir is in a state of equilibrium. After development begins, as fluids are continuously extracted from the reservoir, the formation pressure changes, disrupting the pressure balance relationship [17]. The seepage of fluids in reservoir rocks is also a dynamic process, with state parameters continuously changing with pressure.

(1) Liquid State Equation

Since both oil and water in the reservoir have certain compressibility, during the production process of the oil and gas field, as the reservoir pressure changes, the state of the reservoir fluid will also change. This change is typically described using the liquid compressibility coefficient relationship:

$$\rho = \rho_0 \left[1 + C_l(p - p_0)\right] \qquad (2\text{-}5)$$

Where $\rho$ and $\rho_0$ are the liquid densities at pressures $p$ and $p_0$, respectively, kg/m³; $p_0$ is the atmospheric pressure (or initial pressure), MPa; $C_l$ is the elastic compressibility coefficient of the liquid, MPa⁻¹; $p$ is the reservoir pressure, MPa.

(2) Rock State Equation

Reservoir rocks are classified as porous media and have certain elastic-plastic properties. The relationship between porosity and pore pressure within a certain range of pressure changes can be

expressed as follows:

$$\phi = \phi_0 + C_f(p - p_0) \tag{2-6}$$

Where $\phi$ and $\phi_0$ are the porosities at pressures $p$ and $p_0$, respectively; $p$ is the reservoir pressure, MPa; $p_0$ is the reference pressure, MPa; $C_f$ is the compressibility coefficient of the rock, MPa$^{-1}$.

**2.2 Dual-porous and dual-permeable media model**

The dual-porous and dual-permeable media model was first proposed by Barenblatt in 1960 to study fluid flow in fractured heterogeneous porous media. Its core idea is to divide the reservoir into two overlapping continuous systems: the matrix and the fracture systems. The permeability of the matrix system is several orders of magnitude smaller than its porosity, primarily serving as a fluid storage medium, while the permeability of the fracture system is several orders of magnitude larger than its porosity, mainly facilitating seepage. The flow of fluid within the medium is characterized by the "cross-flow" between these two types of systems [18].

In addressing the issue of fracturing-flooding simulation, it is assumed that there are no open fractures in the initial state of the reservoir, which can be considered a single medium zone [19]. When the water injection well exerts sufficient pressure on the formation, it leads to the opening of fractures, forming a dual-media zone, while the area outside remains the matrix zone. Achieving a dynamic division between this dual-media zone and the matrix zone is challenging. Therefore, we propose an equivalent method, treating the entire reservoir area as a dual-media zone, where the matrix and fracture zones have identical physical property parameters in the initial state, effectively equating them to the matrix zone. When fracturing-flooding begins, different rock stress sensitive variation parameters are set for the fracture and matrix zones to achieve a dynamic description of the dual-media region.

The advantage of the dual-porous and dual-permeable media model lies in its consideration of the physical exchange processes of fluids between the two different systems and the storage effects of fractures in fractured reservoirs, providing significant advantages in describing the dynamic changes of fracture networks in fracturing-flooding reservoirs. The flow patterns among the matrix, fractures, and wellbore in the dual-media model involve supplying fluid from the matrix to the fractures and wellbore, while fluid also flows from the matrix to the fractures and then from the fractures to the wellbore. The fractures serve as the primary fluid flow pathways, while the matrix acts as the main storage space. The fluid flow is characterized as quasi-steady-state flow.

The Darcy seepage model with dual-media characteristics can be expressed by the following equations:

$$\begin{cases} k_f^0 \dfrac{1}{r_D} \dfrac{\partial}{\partial r_D}\left(r_D \dfrac{\partial P_{Df}}{\partial r_D}\right) + \lambda(P_{Dm} - P_{Df}) = \omega_f \dfrac{\partial P_{Df}}{\partial t_D} \\ k_f^0 \dfrac{1}{r_D} \dfrac{\partial}{\partial r_D}\left(r_D \dfrac{\partial P_{Dm}}{\partial r_D}\right) + \lambda(P_{Dm} - P_{Dv}) = \omega_m \dfrac{\partial P_{Dm}}{\partial t_D} \end{cases} \quad (2\text{-}7)$$

Internal boundary conditions:

$$\begin{cases} C_D \dfrac{dP_{wD}}{dt_D} - \left(k_m^0 \dfrac{\partial P_{Dm}}{\partial r_D} + k_f^0 \dfrac{\partial P_{Df}}{\partial r_D}\right)\bigg|_{r_D=1} = 1\,(t_D > 0) \\ P_{wD} = \left(P_{Df} - S\dfrac{\partial P_{Df}}{\partial r_D}\right)\bigg|_{r_D=1} = \left(P_{Dm} - S\dfrac{\partial P_{Dm}}{\partial r_D}\right)\bigg|_{r_D=1} \end{cases} \quad (2\text{-}8)$$

Outer boundary conditions:

$$\lim_{r_D \to \infty} P_{Dm}(r_D, t_D) = \lim_{r_D \to \infty} P_{Df}(r_D, t_D) = 0 \quad (2\text{-}9)$$

Initial conditions:

$$P_{Dj}(r_D, t_D)\big|_{t_D=0} = 0\,(j = m \cdot f)\, 1 \le r_D \le +\infty \quad (2\text{-}10)$$

Where $\omega$ is the storage coefficient, $\lambda$ is the channeling coefficient, and the subscript $D$ is dimensionless.

## 2.3 Anisotropic stress sensitivity

Reservoir stress sensitivity refers to the process in which the original force balance state of solid particles in a porous medium is disrupted due to changes in pore pressure during the extraction of oil and gas reservoirs, leading to the establishment of a new pressure balance state. This phenomenon is also a result of the coupling effects of rock deformation and fluid seepage. After stress sensitivity occurs, the pore space of the reservoir undergoes deformation, altering the seepage effects, which macroscopically manifests as an increase or decrease in permeability, subsequently affecting oil and gas production and leading to changes in output [20]. The dynamic changes in rock properties during the fracturing-flooding process can be represented by curves related to pressure. By incorporating the relationship curve between rock property changes and pressure into the dual-porous and dual-permeable media model, the variation patterns of rock properties during the increase of reservoir pressure can be described, thereby establishing a rock stress sensitive fracturing-flooding simulation method to depict the initiation and dynamic expansion of fractures [21]. When the reservoir pressure changes within a certain limit, the rock properties exhibit partially reversible changes within that limit, which can be represented by a curve. For the matrix, as the pressure increases, both permeability and porosity also increase, but at a relatively gentle rate. When the pressure reaches the fracture pressure point, fractures will form, leading to a rapid increase in permeability, as shown in Fig 2.1 (a). When fractures occur, their permeability is several orders of magnitude greater than that of the porosity [22]. The same porosity versus pressure change curve can be set for both the matrix and fracture zones, as illustrated in Fig 2.1 (b).

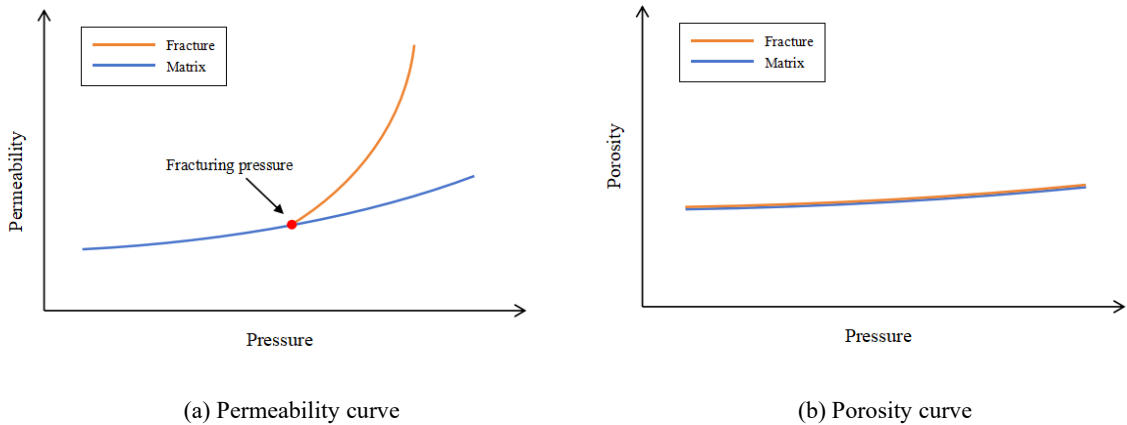

(a) Permeability curve  (b) Porosity curve

Fig 2.1 Schematic diagram of rock stress sensitive curve

Due to the heterogeneity of the reservoir and the influence of geostress, fractures do not expand uniformly in all directions, exhibiting anisotropy [23]. Therefore, the rock stress sensitivity curves in the X and Y directions of the model are different.

**2.4 Implementation of simulation methods**

The above design can be fully implemented in existing commercial numerical simulators such as ECLIPSE and tNavigator, without the need to develop a numerical simulator from scratch, greatly reducing the difficulty of simulation. Taking ECLIPSE as an example:

In the ECLIPSE reservoir simulator, to characterize the rock compaction curve, the compaction curve is typically made into a curve table. It is assumed that the initial reservoir rock is fully compacted, and the initial pore volume and transmissibility multiplier will be set to the values corresponding to the rock compaction curve under the pressure of the grid cell. To demonstrate the different properties of the rock and fractures, two sets of different compaction curve tables are set for the matrix and fracture cells in the model.

The specific steps to establish a model using the ECLIPSE reservoir simulator are as follows:

(1) Set the model as a dual-porosity dual-permeability model by adding the keywords DUALPORO and DUALPERM in the RUNSPEC section of the existing model DATA file.

(2) In the DIMENS keyword, change the dimension of the model in the Z direction to twice that of the original model, where 1-N represents the matrix part; N+1:2N represents the fracture part.

(3) Define the model partition, where the 1-N zone is Region 1, representing the matrix zone; N+1:2N is Region 2, representing the fracture zone.

(4) Assign grid properties to the matrix and fracture zones, giving them the original property values.

(5) Assign relative permeability curves and capillary pressure curves to the matrix and fracture zones, also using the original values.

(6) Assign rock stress sensitivity parameters to the matrix and fracture zones, providing a set of rock stress sensitivity tables for each, which includes a series of characteristic pressure points and the corresponding pore volume multipliers and transmissibility multiplier at each pressure point.

(7) Run the numerical simulation.

## 3 Automatic history matching of model parameters

To achieve a dynamic description of the fracturing-flooding process using the aforementioned numerical simulation method, the accurate setting of rock stress sensitivity parameters is crucial. However, it is often challenging to obtain accurate rock stress sensitivity parameters in field practice, requiring manual adjustments and repeated simulations to achieve an accurate representation of the fracturing-flooding process. To address this, we propose the application of an automatic history matching method to obtain the relevant parameters through the matching of dynamic data from the fracturing-flooding process.

### 3.1 Mathematical model for automatic history matching

(1) Optimization Objectives

The fracturing-flooding development process is divided into three stages: the water injection stage, the soaking stage, and the production stage.

Traditionally, the goal of history matching is to adjust reservoir geological parameters so that the simulated dynamic production indicators, such as bottom hole pressure, water injection rate, liquid production rate, and water cut, are as consistent as possible with the actual monitoring results.

For bottom hole pressure, the objective function for matching can be expressed as:

$$\min \alpha_P = \sum_{w=1}^{M} \sum_{t=1}^{N} \sqrt{\left(P_w^t - P_w^{t'}\right)^2} \qquad (3\text{-}1)$$

Where $\alpha_P$ represents the deviation of the simulated bottom hole pressure for each well, $w$ represents the well number, $M$ represents the total number of wells, $t$ represents the index number of the observation data point, $N$ represents the total number of observation data points, $P_w^t$ is the actual bottom hole pressure at the $t$-th observation point of the $w$-th well, and $P_w^{t'}$ is the simulated bottom hole pressure at the $t$-th observation point of the $w$-th well.

For the water injection rate, the objective function for matching can be expressed as:

$$\min \alpha_Q = \sum_{w=1}^{M} \sum_{t=1}^{N} \sqrt{\left(Q_w^t - Q_w^{t'}\right)^2} \qquad (3\text{-}2)$$

Where $\alpha_Q$ represents the deviation of the simulated water injection rate, $Q_w^t$ is the actual water injection rate at the *t*-th observation point of the *w*-th well, and $Q_w^t$ is the simulated water injection rate at the *t*-th observation point of the *w*-th well.

For oil wells, the produced liquid contains both oil and water phases, necessitating the matching of water cut. The objective function for matching can be expressed as:

$$\min \alpha_f = \sum_{w=1}^{M_o} \sum_{t=1}^{N} \sqrt{\left(f_w^t - f_w^{'t}\right)^2} \tag{3-3}$$

Where $\alpha_f$ represents the deviation of the simulated water cut, $M_o$ represents the total number of oil wells, $f_w^t$ is the actual water cut at the *t*-th observation point of the *w*-th oil well, and $f_w^{'t}$ is the simulated water cut at the *t*-th observation point of the *w*-th oil well.

In fracturing-flooding development, in addition to the above monitoring indicators, it is also necessary to include the matching of the fracture extension range, meaning that the expected fracture extension range obtained from numerical simulation should be as consistent as possible with the actual monitored fracture range. Here, we choose to use the lengths of the maximum and minimum principal stress direction fractures as reference indicators.

The objective function for matching can be expressed as:

$$\min \alpha_{dx} = \sum_{\mu=1}^{M_f} \sum_{t=1}^{N} \sqrt{\left(dx_\mu^t - dx_\mu^{'t}\right)^2} \tag{3-4}$$

$$\min \alpha_{dy} = \sum_{\mu=1}^{M_f} \sum_{t=1}^{N} \sqrt{\left(dy_\mu^t - dy_\mu^{'t}\right)^2} \tag{3-5}$$

Where $\alpha_{dx}$ represents the deviation of the length of the maximum principal stress direction fracture, $\alpha_{dy}$ represents the deviation of the length of the minimum principal stress direction fracture, $M_f$ represents the total number of fractures, $\mu$ represents the fracture number, $dx_\mu^t$ represents the actual length of the $\mu$-th maximum principal stress direction fracture at the *t*-th observation point, $dx_\mu^{'t}$ represents the simulated length of the $\mu$-th maximum principal stress direction fracture at the *t*-th observation point, $dy_\mu^t$ represents the actual length of the $\mu$-th minimum principal stress direction fracture at the *t*-th observation point, and $dy_\mu^{'t}$ represents the simulated length of the $\mu$-th minimum principal stress direction fracture at the *t*-th observation point.

The history matching for fracturing-flooding is a multi-objective optimization problem. To facilitate solving, we choose to convert the multi-objective optimization problem into a single-

objective optimization problem, as expressed in the following transformation formula:

$$\min \alpha = \frac{1}{5}\left(\frac{\alpha_P}{\overline{P}} + \frac{\alpha_Q}{\overline{Q}} + \frac{\alpha_f}{\overline{f}} + \frac{\alpha_{dx}}{\overline{dx}} + \frac{\alpha_{dy}}{\overline{dy}}\right) \qquad (3\text{-}6)$$

Where $\alpha$ represents the total deviation of the simulation parameters, $\overline{P}$ represents the average bottom hole pressure of each observation data point, $\overline{Q}$ represents the average water injection rate of each observation data point, $\overline{f}$ represents the average water cut of each observation data point, $\overline{dx}$ represents the average length of the maximum principal stress direction fractures of each observation data point, and $\overline{dy}$ represents the average length of the minimum principal stress direction fractures of each observation data point.

To address the dimensional imbalance among the various optimization objectives, the deviations are divided by their corresponding average values of the observation data points, thereby achieving dimensional normalization. The average of the dimensionally normalized deviations then represents the total deviation of the simulation parameters.

(2) Optimization variables

Traditionally, the optimization variables for history matching mainly include permeability, porosity, and relative permeability curves. For the history matching of fracturing-flooding, these geological and fluid parameters are already well understood in conventional reservoir numerical simulation studies. Considering the dynamic influencing factors of fracturing-flooding, it is essential to focus on matching rock mechanics-related parameters. Since pressure significantly affects rock mechanics-related parameters, the water compressibility and oil formation volume factor, which are pressure-dependent, also need to be optimized together.

The water compressibility is less influenced by parameters such as formation pressure and can be represented as a scalar value $C_w$. The oil formation volume factor is a curve that varies with pressure and can be represented as a vector $V_o$, satisfying the function relationship that describes its variation with pressure:

$$V_o = g_{V_o}(P) \qquad (3\text{-}7)$$

Where $P$ represents pressure, and $g_{V_o}(P)$ is the function relationship that describes how the oil formation volume factor varies with pressure.

To simplify the optimization variables, this study uses an oil formation volume factor table to describe the curve of the oil formation volume factor as it varies with pressure, as shown in Table 3.1.

Table 3.1 Oil formation volume factor

| Pressure (MPa) | Oil formation volume factor |
|---|---|
| $P_{o1}$ | $V_{o1}$ |
| $P_{o2}$ | $V_{o2}$ |
| …… | …… |
| $P_{on}$ | $V_{on}$ |

In this table, $P_{o1}$, $P_{o2}$,……, $P_{on}$ are the specified $n$ characteristic pressure points, and $V_{o1}$, $V_{o2}$,……, $V_{on}$ are the corresponding oil formation volume factors at those characteristic pressure points.

The rock physical parameters are relatively more complex. In this study, the porosity change is measured using a pore volume multiplier $\lambda$, and the permeability change is measured using a transmissibility multiplier $\Psi$. This includes: the pore volume multiplier curve as it varies with pressure, the X-direction transmissibility multiplier $\Psi_X$ curve as it varies with pressure, the Y-direction transmissibility multiplier $\Psi_Y$ curve as it varies with pressure, and the Z-direction transmissibility multiplier $\Psi_Z$ curve as it varies with pressure, satisfying the function relationship that describes their variation:

Similarly to oil formation volume factor, this study employs a rock stress sensitivity parameter table to describe the relationship of porosity and permeability changes for the matrix and fractures as they vary with pressure. As shown in Table 3.2, two sets of physical property parameter tables need to be established for the matrix and fracture media in the model to separately describe their mechanical properties due to their differing rock mechanics attributes.

Table 3.2 Rock stress sensitivity parameter table

| Parameter | Pressure (MPa) | Pore volume multiplier | X-direction transmissibility multiplier | Y-direction transmissibility multiplier | Z-direction transmissibility multiplier |
|---|---|---|---|---|---|
| Matrix | $P_{M1}$ | $\lambda_{M1}$ | $\Psi_{XM1}$ | $\Psi_{YM1}$ | $\Psi_{ZM1}$ |
|  | $P_{M2}$ | $\lambda_{M2}$ | $\Psi_{XM2}$ | $\Psi_{YM2}$ | $\Psi_{ZM2}$ |
|  | …… |  | …… |  |  |
|  | $P_{Mn}$ | $\lambda_{Mn}$ | $\Psi_{XMn}$ | $\Psi_{YMn}$ | $\Psi_{ZMn}$ |
| Fracture | $P_{F1}$ | $\lambda_{F1}$ | $\Psi_{XF1}$ | $\Psi_{YF1}$ | $\Psi_{ZF1}$ |
|  | $P_{F2}$ | $\lambda_{F2}$ | $\Psi_{XF2}$ | $\Psi_{YF2}$ | $\Psi_{ZF2}$ |
|  | …… |  | …… |  |  |
|  | $P_{Fn}$ | $\lambda_{Fn}$ | $\Psi_{XFn}$ | $\Psi_{YFn}$ | $\Psi_{ZFn}$ |

Where $P_{M1}$, $P_{M2}$,……, $P_{Mn}$ are the specified n$n$ characteristic pressure points for the matrix; $\lambda_{M1}, \lambda_{M2}$,……,$\lambda_{Mn}$ are the pore volume multiplier of the matrix at the corresponding characteristic pressure points; $\Psi_{XM1}$, $\Psi_{XM2}$, ……, $\Psi_{XMn}$ are the X-direction transmissibility multipliers of the matrix at the corresponding characteristic pressure points; $\Psi_{YM1}$, $\Psi_{YM2}$, ……, $\Psi_{YMn}$ are the Y-direction transmissibility multipliers of the matrix at the corresponding characteristic pressure points; and $\Psi_{ZM1}$, $\Psi_{ZM2}$, ……, $\Psi_{ZMn}$ are the Z-direction transmissibility multiplier of the matrix at the corresponding characteristic pressure points. $P_{F1}$, $P_{F2}$, ……, $P_{Fn}$ are the specified n$n$ characteristic pressure points for the fractures; $\lambda_{F1}$, $\lambda_{F2}$, ……, $\lambda_{Fn}$ are the pore volume multiplier of the fractures at the corresponding characteristic pressure points; $\Psi_{XF1}$, $\Psi_{XF2}$, ……, $\Psi_{XFn}$ are the X-direction transmissibility multipliers of the fractures at the corresponding characteristic pressure points; $\Psi_{YF1}$, $\Psi_{YF2}$, ……, $\Psi_{YFn}$ are the Y-direction transmissibility multipliers of the fractures at the corresponding characteristic pressure points; and $\Psi_{ZF1}$, $\Psi_{ZF2}$, ……, $\Psi_{ZFn}$ are the Z-direction transmissibility multipliers of the fractures at the corresponding characteristic pressure points.

Thus, when $n$ characteristic pressure points are given, there are $2n$ optimization variables in the oil formation volume factor table, $10n$ optimization variables in the rock stress sensitivity table, and a separate optimization variable $C_w$. The total number of optimization variables $T$ is:

$$T = 12n + 1 \tag{3-10}$$

The process of implementing fracturing-flooding consists of three stages, and to accurately match the rock physical property change curve, at least three characteristic pressure points should be set, making $T$ no less than 37. It can be seen that if matching is done directly, there are still many variables to optimize, and the relationships between the parameters are complex. To facilitate calculations, further simplification of the optimization variables is necessary.

For the oil formation volume factor table, a series of initial values that conform to the experimental rules can be set for all oil formation volume factors based on field test understanding, and an overall coefficient multiplier $K_{Vo}$ can be used to represent their magnitude. For the rock stress sensitivity table, taking the pore volume multiplier as an example, three characteristic pressure points are set to approximate the variation relationship. As shown in Fig 3.1, pressure points $P_{min}$ (lower limit pressure point, known), $P_{max}$ (upper limit pressure point, known), and $P_\sigma$ (intermediate pressure point) are established, where $P_{min}$ and $P_{max}$ can be specified based on field understanding, while the value of $P_\sigma$ needs to be matched. If the pore volume multiplier corresponding to these three characteristic pressure points are matched directly, inequality constraints need to be applied ($\lambda_{min} < \lambda_\sigma < \lambda_{max}$). To simplify the constraint conditions, the optimization variables are transformed into the increments $\Delta\lambda_1$ and $\Delta\lambda_2$ between the pore volume multiplier and $\lambda_{min}$ (the lower limit pore volume multiplier). Similarly, this approach is applied to the transmissibility multipliers, which can achieve the goal of reducing the total number of optimization variables and facilitate the description of the size relationships among the optimization variables.

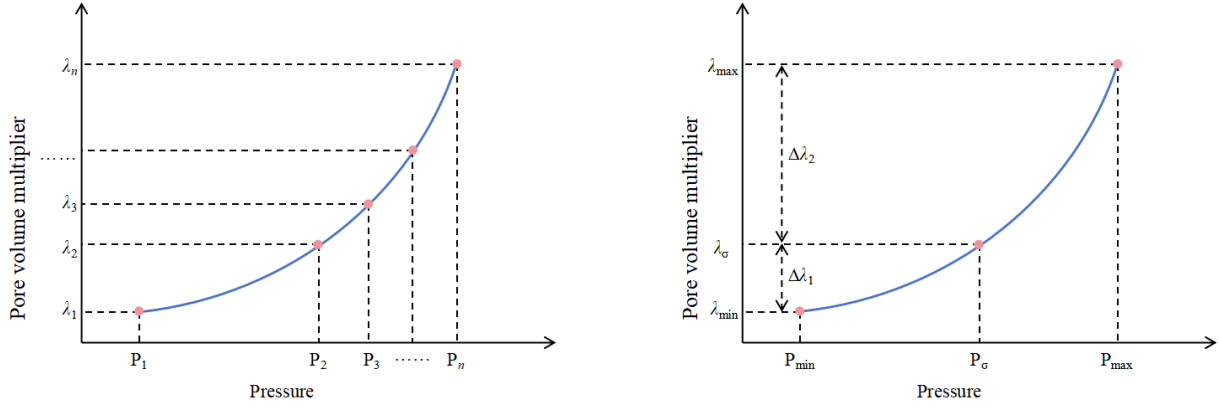

(a) Before simplification  (b) After simplification

Fig 3.1 Simplification of pore volume multiplier

Considering that the physical properties of the matrix and fractures in the model change consistently before the fracture pressure, the optimization parameters for the matrix and fractures at the first two characteristic pressure points can be integrated. This way, only the characteristic points after the fracture pressure point $P_B$ have inconsistent optimization parameters. Due to the influence of geostress anisotropy, the permeability in the X, Y, and Z directions will vary. Based on the X-direction transmissibility multiplier, this study calculates the transmissibility in other directions using an anisotropy coefficient $K_{XY}$.

Finally, the following representative parameters are defined: water compressibility coefficient $C_w$, oil formation volume factor $K_{Vo}$, fracture pressure point $P_B$, lower limit rock pore volume multiplier $\lambda_{Mmin}$, pore volume multiplier increments $\Delta\lambda_1$ and $\Delta\lambda_2$, lower limit rock transmissibility multiplier $\Psi_{XMmin}$, transmissibility multiplier increments $\Delta\Psi_{XM1}$ and $\Delta\Psi_{XM2}$, maximum stress sensitive direction transmissibility multiplier in the fracture zone $\Psi_{XFmax}$, and anisotropy coefficient $K_{XY}$. As shown in Table 3.3 and Table 3.4, based on these 11 representative parameters, a stress sensitivity table can be automatically generated.

Table 3.3 Simplified Oil formation volume factor

| Pressure (MPa) | Oil formation volume factor |
|---|---|
| $P_{O1}$ | $V_{O1} \cdot K_{Vo}$ |
| $P_{O2}$ | $V_{O2} \cdot K_{Vo}$ |
| …… | …… |
| $P_{On}$ | $V_{On} \cdot K_{Vo}$ |

Table 3.4 Simplified Rock stress sensitivity parameter table

| Parameter | Pressure | Pore volume | X-direction | Y-direction | Z-direction |
|---|---|---|---|---|---|

|  | (MPa) | multiplier | transmissibility multiplier | transmissibility multiplier | transmissibility multiplier |
|---|---|---|---|---|---|
| Matrix | $P_{min}$ | $\lambda_{Mmin}$ | $\Psi_{XMmin}$ | $\Psi_{XMmin}$ | $\Psi_{XMmin}$ |
|  | $P_B$ | $\lambda_{Mmin}+\Delta\lambda_1$ | $\Psi_{XMmin}+\Delta\Psi_{XM1}$ | $\Psi_{XMmin}+\Delta\Psi_{XM1}$ | $\Psi_{XMmin}+\Delta\Psi_{XM1}$ |
|  | $P_{max}$ | $\lambda_{Mmin}+\Delta\lambda_1+\Delta\lambda_2$ | $\Psi_{XMmin}+\Delta\Psi_{XM1}+\Delta\Psi_{XM2}$ | $\Psi_{XMmin}+\Delta\Psi_{XM1}+\Delta\Psi_{XM2}$ | $\Psi_{XMmin}+\Delta\Psi_{XM1}+\Delta\Psi_{XM2}$ |
| Fracture | $P_{min}$ | $\lambda_{Mmin}$ | $\Psi_{XMmin}$ | $\Psi_{XMmin}$ | $\Psi_{XMmin}$ |
|  | $P_B$ | $\lambda_{Mmin}+\Delta\lambda_1$ | $\Psi_{XMmin}+\Delta\Psi_{XM1}$ | $\Psi_{XMmin}+\Delta\Psi_{XM1}$ | $\Psi_{XMmin}+\Delta\Psi_{XM1}$ |
|  | $P_{max}$ | $\lambda_{Mmin}+\Delta\lambda_1+\Delta\lambda_2$ | $\Psi_{XFmax}$ | $K_{XY}\cdot\Psi_{XFmax}$ | $K_{XY}\cdot\Psi_{XFmax}$ |

(3) Constraints

In optimization problems, it is necessary to set constraints on the variables to ensure that the solution space remains within a reasonable range. As scalar values, $C_w$ and $K_{Vo}$ need to have upper and lower bounds established. For the variables in the oil formation volume factor table and the rock stress sensitivity table, not only must the variable ranges be defined, but the relationships between the sizes of each variable must also be constrained. The specific relationships are as follows:

$$\begin{cases} P_{M1}<P_{M2}<\cdots\cdots<P_{Mn} \\ \lambda_{M1}<\lambda_{M2}<\cdots\cdots<\lambda_{Mn} \\ \Psi_{XM1}<\Psi_{XM2}<\cdots\cdots<\Psi_{XMn} \\ \Psi_{YM1}<\Psi_{YM2}<\cdots\cdots<\Psi_{YMn} \\ \Psi_{ZM1}<\Psi_{ZM2}<\cdots\cdots<\Psi_{ZMn} \\ P_{F1}<P_{F2}<\cdots\cdots<P_{Fn} \\ \lambda_{F1}<\lambda_{F2}<\cdots\cdots<\lambda_{Fn} \\ \Psi_{XF1}<\Psi_{XF2}<\cdots\cdots<\Psi_{XFn} \\ \Psi_{YF1}<\Psi_{YF2}<\cdots\cdots<\Psi_{YFn} \\ \Psi_{ZF1}<\Psi_{ZF2}<\cdots\cdots<\Psi_{ZFn} \end{cases} \quad (3\text{-}11)$$

Describing the above inequality relationships in the matching algorithm is quite difficult and complex to solve. Therefore, the aforementioned method simplifies the number of variables while also transforming the relationships between the sizes of the variables. It solves for the minimum volume multiplier and transmissibility multiplier corresponding to $P_{min}$, and sets a series of increments ($\Delta\lambda_1$, $\Delta\lambda_2$, $\Delta\Psi_{XM1}$, $\Delta\Psi_{XM2}$) that are all greater than zero as the differences between adjacent variables. After simplification, a few variables that only require upper and lower bounds replace the complex inequality constraint relationships.

The specific value ranges for each variable need to be determined based on the actual conditions of the target reservoir. According to existing knowledge, the reference values for the constraints of each optimization variable are as follows:

$$\begin{cases} C_w \in [1\times 10^{-6}, 1\times 10^{-4}] \\ K_{V_o} \in [0.8, 1.5] \\ P_B \in [P_{min}, P_{max}] \\ \lambda_{Mmin} \in [0.9, 0.99] \\ \Delta\lambda_1 \in [0.001, 0.05] \\ \Delta\lambda_2 \in [0.001, 0.05] \\ \Psi_{XMmin} \in [0.9, 0.99] \\ \Delta\Psi_{XM1} \in [0.001, 0.05] \\ \Delta\Psi_{XM2} \in [0.001, 0.05] \\ \Delta\Psi_{XFmax} \in [100, 2000] \\ K_{XY} \in [0.1, 0.6] \end{cases} \tag{3-12}$$

**3.2 Solving method of model**

In this study, automatic history matching is designed as an iterative solving problem with multiple optimization variables and a single optimization objective, subject only to upper and lower bound constraints. Although there is only one optimization objective, the solution space exhibits high complexity and nonlinearity, making the use of a gradient-free optimization algorithm more reasonable. The Covariance Matrix Adaptation Evolution Strategy (CMA-ES) is a fast and efficient gradient-free optimization algorithm developed based on evolutionary strategies, characterized by strong search capability and robustness, and is suitable for solving optimization variables constrained by upper and lower bounds [24]. Given its characteristics and advantages, this study selects CMA-ES to solve the mathematical model of automatic history matching.

The CMA-ES algorithm is a population-based random search algorithm, where the individuals in the CMA-ES population are distributed according to a specific probability distribution, and the main adjustments during the iterative optimization process are made to this probability distribution.

In the iteration step $k$, CMA-ES first samples $\gamma$ individuals to form a population according to the following formula:

$$\mathbf{x}_i^k = N\left(m^k, \left(\sigma^k\right)^2 C^k\right), \text{ for } i = 1, \cdots, \gamma \tag{3-13}$$

Where N (…,…) is a random vector from a multivariate normal distribution; $m^k$ is the mean vector; $C^k$ is the covariance matrix; and $\sigma$ is the step size factor.

The mean vector $m^k$ represents the current optimal solution; the covariance matrix $C^k$ is a symmetric positive definite matrix used to describe the geometric characteristics of the distribution; the step size factor $\sigma$ is used to globally scale the covariance matrix $C^k$ to achieve rapid convergence and avoid premature convergence. During the iteration process, CMA-ES requires continuous updates of these three parameters.

The mean vector $m^k$ is obtained by calculating the weighted average of the $\mu$ individuals with the smallest objective function values, with the calculation formula as follows:

$$m^{k+1} = \sum_{i=1}^{\mu} \omega_i x_{1:\gamma}^{k} \tag{3-14}$$

The default weights are:

$$\omega_i = \frac{\ln(\mu+1) - \ln(i)}{\mu \ln(\mu+1) - \ln(\mu!)}, \text{ for } i = 1, \cdots, \mu \tag{3-15}$$

In general, $\mu$ is set to equal $\mu = \lfloor \lambda/2 \rfloor$, where $\lfloor \ \rfloor$ is the floor function, and ωi are all positive values that sum to 1.

Subsequently, the covariance matrix is updated as follows:

$$C^{k+1} = (1 - c_{cov})C^k + \frac{c_{cov}}{\mu_{cov}} p_c^{k+1} p_c^{(k+1)T} + c_{cov}\left(1 - \frac{1}{\mu_{cov}}\right) \times \sum_{i=1}^{\mu} \frac{\omega_i}{\sigma^{(k)2}} \left(x_{i:1}^{k+1} - m^k\right)\left(x_{i:1}^{k+1} - m^k\right)^T \tag{3-16}$$

Where $p_c^k$ becomes the evolution path, representing the optimization direction.

$p_c^k$ also needs to be updated at each iteration step, with the update formula as follows:

$$p_c^{k+1} = (1 - c_c)p_c^k + \sqrt{c_c(2 - c_c)\mu_{eff}} \frac{m^{k+1} - m^k}{\sigma^k} \tag{3-17}$$

Where $c_c$ is a constant between 0 and 1; $\mu_{eff} = 1/\sum_{i=1}^{\mu}\omega_i^2$ is used to describe recombination, and the new search direction $p_c^{k+1}$ is determined by both the original search direction $p_c^{k+1}$ and the descent direction $\frac{m^{k+1} - m^k}{\sigma^k}$.

The step size factor $\sigma^{k+1}$ is determined by the following formula:

$$\sigma^{k+1} = \sigma^k \exp\left[\frac{c_\sigma}{d_\sigma}\left(\frac{p_\sigma^{k+1}}{E\|N(0,I)\|} - 1\right)\right] \tag{3-18}$$

Where $p_\sigma^{k+1}$ is the conjugate search path, determined as follows:

$$p_\sigma^{k+1} = (1 - c_\sigma)p_\sigma^k + \sqrt{c_\sigma(2 - c_\sigma)\mu_{eff}} \sigma^{k-1} C^{k-\frac{1}{2}} \left(m^{k+1} - m^k\right) \tag{3-19}$$

By adaptively scaling the covariance matrix with the step size factor, good convergence speed can be achieved.

**3.3 Overall process**

The overall solving process of automatic history matching based on the CMA-ES algorithm is shown in Fig 3.2. First, initial values for each optimization variable that meet the constraint conditions are set and imported into the corresponding positions in the model file, and the model is run. A pre-defined program is used to extract and process the data from the result file, and the matching situation of the result data is evaluated. The CMA-ES optimization algorithm is then called to update the optimization variables based on the matching situation, and the updated parameters are imported into the model for running. This process is iterated until the optimal parameter combination is found.

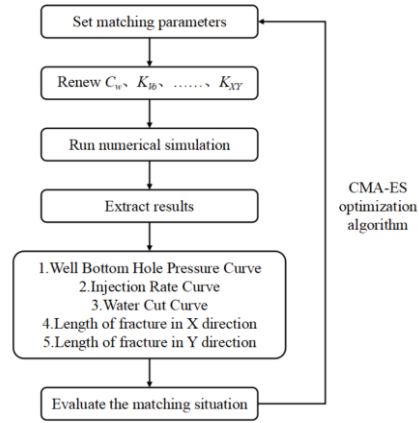

Fig 3.2 Solving process of automatic history matching

# 4 Conclusion

Based on the theories of seepage mechanics and rock mechanics, this paper establishes a dual-porous and dual-permeable media model that is stress sensitive to rock anisotropy. By employing an automatic history matching method based on the CMA-ES optimization algorithm, model parameters that align with actual production conditions are derived, resulting in a high-precision model suitable for fracturing-flooding production. In conjunction with the implementation status of fracturing-flooding in the target block, the dynamic changes of fractures during the fracturing-flooding process are studied, leading to the following conclusions:

(1) This paper introduces a Darcy flow model with dual-porous and dual-permeable characteristics based on the theory of seepage mechanics. It employs two sets of rock stress-sensitive parameter tables to describe the physical property changes of the matrix and fractures in the reservoir during the fracturing-flooding process, setting different parameters in the X and Y directions to characterize the anisotropic features of low-permeability reservoirs. This model can accurately describe the dynamic changes of fractures in the reservoir during the fracturing-flooding process, including the geometric shape of fracture network extension and the expansion and closure of the fracture network as production evolves.

(2) An automatic history matching method based on the CMA-ES algorithm has been designed. This matching method uses the bottom hole pressure curves, water injection rate curves, water cut curves, and fracture extension ranges of each production well during the fracturing-flooding process as optimization objectives, while reservoir physical property parameters, such as pore volume multipliers, are used as optimization variables. The CMA-ES algorithm is employed to solve for the optimal parameter combinations. This method can improve modeling efficiency, enhance matching accuracy, and reduce labor costs.